\documentclass[11pt]{article}
\usepackage{amsfonts,latexsym,amsmath,amscd}
\newcommand{\beq}{\begin{equation}}
\newcommand{\eeq}{\end{equation}}
\newcommand{\beqa}{\begin{eqnarray}}
\newcommand{\eeqa}{\end{eqnarray}}
\newcommand{\beaa}{\begin{eqnarray*}}
\newcommand{\ben}{\begin{eqnarray*}}
\newcommand{\eaa}{\end{eqnarray*}}
\newcommand{\een}{\end{eqnarray*}}

\newcommand{\bal}{\begin{align}}
\newcommand{\eal}{\end{align}}

\newcommand \nc {\newcommand}
\nc \proof {\noindent {\em{Proof.\/ }}}
\nc \qed {\hfill $\Box$}
\newtheorem{theorem}{Theorem}[section]
\newtheorem{lemma}[theorem]{Lemma}
\newtheorem{proposition}[theorem]{Proposition}
\newtheorem{corollary}[theorem]{Corollary}
\newtheorem{definition}[theorem]{Definition}
\newtheorem{example}[theorem]{Example}
\newtheorem{remark}[theorem]{Remark}
\newtheorem{conjecture}[theorem]{Conjecture}
\newtheorem{question}[theorem]{Question}
\nc \bth[1] {\begin{theorem}\label{t#1} }
\nc \ble[1] {\begin{lemma}\label{l#1} }
\nc \bpr[1] {\begin{proposition}\label{p#1} }
\nc \bco[1] {\begin{corollary}\label{c#1} }
\nc \bde[1] {\begin{definition}\label{d#1}\rm }
\nc \bex[1] {\begin{example}\label{e#1}\rm }
\nc \bre[1] {\begin{remark}\label{r#1}\rm }
\nc \bcon[1] {\begin{conjecture}\label{con#1}\rm }
\nc \bque[1] {\begin{question}\label{que#1}\rm }
\nc {\eth} {\end{theorem}}
\nc {\ele} {\end{lemma}}
\nc {\epr} {\end{proposition}}
\nc {\eco} {\end{corollary}}
\nc {\ede} {\end{definition}}
\nc {\eex} {\end{example}}
\nc {\ere} {\end{remark}}
\nc {\econ} {\end{conjecture}}
\nc {\eque} {\end{question}}
\nc \thref[1]{Theorem \ref{t#1}}
\nc \leref[1]{Lemma \ref{l#1}}
\nc \prref[1]{Proposition \ref{p#1}}
\nc \coref[1]{Corollary \ref{c#1}}
\nc \deref[1]{Definition \ref{d#1}}
\nc \exref[1]{Example \ref{e#1}}
\nc \reref[1]{Remark \ref{r#1}}

\def \( {{\left(}}
\def \) {{\right)}}

\def \[  {{\left[}}
\def \]  {{\right]}}

\def\tensor{\otimes}

\begin{document}
\title{A Remark on Nonlinear Dirac Equations}
\author{Changyou Wang\thanks{Partially supported by NSF grant 0601162}\\
Department of Mathematics\\
University of Kentucky\\
Lexington, KY 40506, USA\\
{\it cywang@ms.uky.edu}}
\date{}
\maketitle

\begin{abstract}
For a $n$-dimensional spin manifold $M$ with a fixed spin structure and a spinor bundle $\Sigma M$,
we prove an $\epsilon$-regularity theorem for weak solutions to the nonlinear Dirac equation
of cubic nonlinearity. This, in particular, answers a regularity question raised by Chen-Jost-Wang
\cite{chen-jost-wang2} when $n=2$.
\end{abstract}

\section{Introduction}

Linear Dirac type equations, including the Cauchy-Riemann equation in dimension two,
are the most fundamental first order system of elliptic equations. During the course
to study {\it Dirac-harmonic maps with curvature term} from a Riemann surface into a
Riemannian manifold,  Chen-Jost-Wang \cite{chen-jost-wang1, chen-jost-wang2} introduced  the nonlinear
Dirac equation with cubic nonlinearity:
\beq{}\label{dirac0}
\partial\hspace{-0.25cm}/\hspace{+0.05cm}\psi^i=\sum_{j,k,l=1}^N H_{jkl}^i\langle \psi^j, \psi^k\rangle \psi^l, \ 1\le i\le N.
\eeq
In dimension two, an interesting feature of this nonlinear Dirac equation is
that it is conformally invariant and has critical nonlinearity, where the classical
methods fail to apply. Thus it is an interesting question to study the regularity of
weak solutions of (\ref{dirac}). The aim of this short note is to provide an elementary
proof of a general regularity criterion for (\ref{dirac0}).

In order to describe the results, we briefly review some background materials on spin manifolds.
The interested readers can consult with Lawson-Michelsohn \cite{lawson-michelsohn}, Chen-Jost-Li-Wang
\cite{chen-jost-li-wang1, chen-jost-li-wang2} for more details.
For $n\ge 2$, let $(M,g)$ be a spin manifold with a given spin structure and an associated
spinor bundle $\Sigma$.  Let $\langle\cdot,\cdot\rangle$ be a Hermitian metric on $\Sigma$ and
$\nabla$ be the Levi-Civita connection on $\Sigma$ compatible with both $\langle\cdot,\cdot\rangle$
and $g$. The Dirac operator on $M$ is defined by
$\partial\hspace{-0.25cm}/\hspace{+0.05cm}=e_\alpha\circ\nabla_{e_\alpha}$,
where $\{e_\alpha\}_{\alpha=1}^n$ is a local orthonormal frame on $M$, and $\circ: TM\tensor_{\mathbb C}\Sigma\to\Sigma$
is the Clifford multiplication.

Now let's write (\ref{dirac}) into the form
\beq{}\label{dirac}
\partial\hspace{-0.25cm}/\hspace{+0.05cm}\psi= H_{jkl}\langle \psi^j, \psi^k\rangle \psi^l,
\eeq
where $\psi=(\psi^1,\cdots,\psi^N)\in \left(\Gamma\Sigma\right)^N$, $N\ge 1$, 
$H_{jkl}=\left(H_{jkl}^1,\cdots, H_{jkl}^N\right)\in C^\infty(M,\mathbb R^N)$.

We refer the readers to \cite{chen-jost-wang2} \S1, where the authors discussed two interesting
examples in which (\ref{dirac}) arises naturally. The first example is the Dirac-harmonic map
$(\phi,\psi)$ associated with the Dirac-harmonic energy functional with curvature term,
a nonlinear $\sigma$-model in the superstring theory, in which the nonlinear Dirac equation
for $\psi$ reduces to (\ref{dirac}) when $\phi$ is a constant map. The second example is
the Weierstrass representation formula for minimal surfaces $X$ immersed in $\mathbb R^3$
by holomorphic $1$-forms and meromorphic functions, in which an equation of the form (\ref{dirac})
appears naturally.

It turns out that the underlying function space for the equation (\ref{dirac}) is $L^4(M)$. As pointed
out by \cite{chen-jost-wang2} that any weak solution $\psi$ of (\ref{dirac}) is smooth provide
$\psi\in L^p(M)$ for some $p>4$.  In \cite{chen-jost-wang2}, the authors proved three interesting
analytic properties of (\ref{dirac}) for $n=2$: (i) the gradient estimate for {\it smooth} solutions $\psi$
of (\ref{dirac}) under the smallness condition of $L^4$-norm of $\psi$, (ii) the isolated singularity
removable theorem, and (iii) the energy identity theorem for sequentially weak convergent smooth solutions
of (\ref{dirac}).  At the end of \S1 in \cite{chen-jost-wang2}, the authors asked

\bcon{}\label{full_reg}
{\em For $n=2$, any weak solution $\psi\in L^4(M)$ of (\ref{dirac}) is smooth.}
\econ

In this short note, we answer Conjecture \ref{full_reg} affirmatively. In fact, we prove a general regularity
theorem for weak solutions of (\ref{dirac}) in any dimensions. The ideas is based on an application of the
estimate of Reisz potentials between Morrey spaces, due to Adams \cite{adams}. Similar techniques have
been employed in the proof of higher order regularity of Dirac-harmonic maps by Wang-Xu \cite{wang-xu}.
The proof turns out to be very elementary, and may be applicable to other similar problems.

Before stating our results, let's first recall
the definition of weak solutions of (\ref{dirac}).
\bde{}\label{weak_sol} A section $\psi\in L^4((\Gamma\Sigma)^N)$ is a weak solution of (\ref{dirac}) if
\beq{}
\label{weak_sol1}
\int_M \langle\psi, \partial\hspace{-0.25cm}/\hspace{+0.05cm}\eta\rangle
=\int_M H_{jkl}\left\langle \psi^j, \psi^k\right\rangle\left\langle \psi^l, \eta\right\rangle
\eeq
holds for any smooth section $\eta\in C^\infty\left((\Gamma\Sigma)^N\right)$.
\ede

Denote by ${\it i}_M>0$ the injectivity radius of $M$. For $0<r<{\it i}_M$ and $x\in M$, denote by $B_r(x)$
the geodesic ball in $M$ with center $x$ and radius $r$. Now we state our theorems.
\bth{}\label{epsilon_reg} For any $n\ge 2$, there exists $\epsilon_0>0$ depending on $n$ such that
if $\psi\in L^4((\Gamma\Sigma)^N)$ is a weak solution of the Dirac equation
(\ref{dirac}) and satisfies, for some $x_0\in M$ and $0<r_0\le \frac12{\it i}_M$,
\beq{} \label{small_cond}
\sup_{x\in B_{r_0}(x_0), \ 0<r\le r_0}\left\{\frac{1}{r^{n-2}}\int_{B_r(x)} |\psi|^4 \right\}\le\epsilon_0^4,
\eeq
then $\psi\in C^\infty(B_{\frac{r_0}2}(x_0))$.
\eth

Note that by H\"older inequality, we have for $n\ge 2$,
$$\frac{1}{r^{n-2}}\int_{B_r(x)} |\psi|^4 \le \left(\int_{B_r(x)} |\psi|^{2n}\right)^{\frac{2}{n}}.$$
Thus, as an immediate consequence of Theorem \ref{epsilon_reg}, we obtain
\bco{} \label{l-2n-reg} For $n\ge 2$, if $\psi\in L^{2n}((\Gamma\Sigma)^N)$ is a weak solution of
the Dirac equation (\ref{dirac}), then $\psi\in C^\infty((\Gamma\Sigma)^N)$.
\eco
It is clear that when $n=2$, Corollary \ref{l-2n-reg} implies Conjecture \ref{full_reg}.

\section{Proof of Theorem \ref{epsilon_reg}}

This section is devoted to the proof of Theorem \ref{epsilon_reg}. Since the regularity is a local property,
we assume, for simplicity of presentation, that for $x_0\in M$, the geodesic ball $B_{{\it i}_M}(x_0)\subset M$
with the metric $g$ is identified by $(B_2,g_0)$. Here $B_2\subset\mathbb R^n$ is the ball with center $0$
and radius $2$, and $g_0$ is the Euclidean metric on $\mathbb R^n$. We also assume that
$\Sigma\big|_{B_2}=B_2\times \mathbb C^L$, with $L={\rm{rank}}_{\mathbb C}\Sigma$.

Let's also recall the definition of Morrey spaces.
\bde{}\label{morrey} For $1\le p\le n$, $0<\lambda\le n$, and a domain $U\subseteq\mathbb R^n$, the
Morrey space $M^{p,\lambda}(U)$ is defined by
$$
M^{p,\lambda}(U)
:=\left\{f\in L^p_{\hbox{loc}}(U): \|f\|_{M^{p,\lambda}(U)}<+\infty\right\},$$
where
$$\left\Vert f\right\Vert _{M^{p,\lambda}(U)}^p=\sup\left\{r^{\lambda-n}\int_{B_r}|f|^p:\ B_r\subseteq U\right\}.$$
It is easy to see that for $1\leq p\leq n$,  $M^{p,\lambda}(U)\subset L^p(U)$,
$M^{p,n}(U)=L^p(U)$ and  $M^{p,p}(U)$ behaves like $L^{n}(U)$ from the view of scalings.
\ede

It is clear that the condition (\ref{small_cond}) in Theorem \ref{epsilon_reg} is equivalent to
$$\left\|\psi\right\|_{M^{4,2}(B_{r_0}(x_0))}\le \epsilon_0.$$
Thus Theorem \ref{epsilon_reg} follows from the following Lemma.
\ble{}\label{small_integ} For any $4<p<+\infty$ and $n\ge 2$,
there exists $\epsilon_0>0$ depending only on $p$ and $n$
such that if $\psi\in M^{4,2}(B_1)$ is a weak solution of (\ref{dirac}) and
$$\|\psi\|_{M^{4,2}(B_1)}\le\epsilon_0,$$
then $\psi\in L^p(B_\frac1{16},\mathbb C^{NL})$. Furthermore, $\psi\in C^\infty(B_\frac1{16}, \mathbb C^{NL})$ and
the estimate 
\beq{}\label{prior-est}
\left\|\nabla^l\psi\right\|_{C^0(B_{\frac{1}{16}})}\le C(\epsilon_0, l), \ \forall l\ge 1
\eeq
holds
\ele
\proof Applying $\partial\hspace{-0.25cm}/\hspace{+0.1cm}$ to (\ref{dirac}), we have, for $1\le i\le N$,
\beq{}\label{dirac1}
\partial\hspace{-0.25cm}/\hspace{+0.05cm}^2\psi^i
=\partial\hspace{-0.25cm}/\hspace{+0.05cm}\left(H_{jkl}^i\langle\psi^j, \psi^k\rangle \psi^l\right)
\eeq
in the sense of distributions. By Lichnerowitz's formula (cf. \cite{lawson-michelsohn}), we have
$$-\Delta\psi^i=\partial\hspace{-0.25cm}/\hspace{+0.05cm}^2\psi^i.$$
Hence we have
\beq{}\label{dirac2}
-\Delta\psi^i
=\partial\hspace{-0.25cm}/\hspace{+0.05cm}\left(H_{jkl}\langle\psi^j, \psi^k\rangle \psi^l\right)
\eeq
in the sense of distributions.

For $m=1,2$, let $\eta_m\in C_0^\infty(B_1)$ be such that $0\le\eta_m\le 1$, $\eta_m\equiv 1$ on
$B_{2^{1-2m}}$. For $1\le i\le N$, define $f^i_m:\mathbb R^n\to\mathbb C^L$ by letting
\beq{}\label{auxil1}
f^i_m(x)=\int_{\mathbb R^n}\frac{\partial G(x,y)}{\partial y_\alpha}
\frac{\partial}{\partial y_\alpha}\circ \left(\eta_m^3H_{jkl}\langle\psi^j, \psi^k\rangle \psi^l\right)
(y)\,dy,\eeq
where $G(x,y)$ is the fundamental solution of $\Delta$ on $\mathbb R^n$.  For $1\le i\le N$, define $g_m^i: B_1\to\mathbb C$ by
letting
\beq{}\label{auxil2} \psi^i=f_m^i+g_m^i.
\eeq
Direct calculations imply that for $m=1,2$ and $1\le i\le N$,
\begin{eqnarray}\label{dirac3}
-\Delta f^i_m&=&\partial\hspace{-0.25cm}/\hspace{+0.05cm}
\left(\eta_m^3H_{jkl}\langle\psi^j, \psi^k\rangle \psi^l\right)\nonumber\\
&=&\partial\hspace{-0.25cm}/\hspace{+0.05cm}\left(
H_{jkl}\langle\psi^j, \psi^k\rangle \psi^l\right) \ \hbox{ in } B_{2^{1-2m}}.
\end{eqnarray}
This and (\ref{dirac2}) imply
\beq{}\label{dirac4}
\Delta g^i_m=0 \ \hbox{ in }B_{2^{1-2m}}.
\eeq
It follows from (\ref{auxil1}) that for $m=1,2$ and $1\le i\le N$,
\beq{}\label{one_reisz}
\left|f^i_m\right|(x)\le C\int_{\mathbb R^n}\left|x-y\right|^{1-n}\left(\eta_m(y)|\psi(y)|\right)^3\,dy
=CI_1(\eta_m^3|\psi|^3)(x),
\eeq
where
$$I_1(f)(x)=\int_{\mathbb R^n}\left|x-y\right|^{1-n}|f(y)|\,dy,\ \ f:\mathbb R^n\to \mathbb R,$$
is the Riesz potential of order one. Let's recall Adams' inequality on Morrey spaces (cf. \cite{adams}):
\beq{}\label{adams}
\left\|I_1(f)\right\|_{M^{\frac{\lambda q}{\lambda-q},\lambda}(\mathbb R^n)}
\le C\|f\|_{M^{q,\lambda}(\mathbb R^n)},  \ \forall 1\le q<\lambda\le n.
\eeq

\noindent{\it Step} 1 ($m=1$).
Since $(\eta_1|\psi|)^3\in M^{\frac43,2}(\mathbb R^n)$, (\ref{adams}) implies that
for $1\le i\le N$,
\begin{eqnarray}\label{morrey_est}
\|f^i_1\|_{M^{4,2}(\mathbb R^n)}
&\le & C\|\eta_1^3|\psi|^3\|_{M^{\frac43,2}(\mathbb R^n)}
=C\|\eta_1|\psi|\|_{M^{4,2}(\mathbb R^n)}^3\nonumber\\
&\le& C\|\psi\|_{M^{4,2}(B_1)}^3\le C\epsilon_0^2\|\psi\|_{M^{4,2}(B_1)}.
\end{eqnarray}
On the other hand, by the standard estimate for harmonic functions, we have that for any
$\theta\in (0,\frac14)$ and $x_0\in B_{\frac14}$
\beq{}\label{morrey_est1}
\|g^i_1\|_{M^{4,2}(B_\theta(x_0))}\le C\theta^{\frac12} \|g^i_1\|_{M^{4,2}(B_{\frac12})},
\ \forall 1\le i \le N.
\eeq
Putting (\ref{morrey_est}) and (\ref{morrey_est1}) into (\ref{auxil2}) yields that for $1\le i\le N$,
\begin{eqnarray}
\|\psi^i\|_{M^{4,2}(B_\theta(x_0))}
&\le & C\theta^{\frac12}\|g^i_1\|_{M^{4,2}(B_{\frac12})}+C\epsilon_0^2\|\psi\|_{M^{4,2}(B_1)}\nonumber\\
&\le & C\theta^{\frac12}\left[\|\psi^i\|_{M^{4,2}(B_{\frac12})}
+\|f^i_1\|_{M^{4,2}(B_{\frac12})}\right]
+C\epsilon_0^2\|\psi\|_{M^{4,2}(B_1)}\nonumber\\
&\le & C\left(\epsilon_0^2+\theta^{\frac12}\right)\|\psi\|_{M^{4,2}(B_1)}.
\end{eqnarray}
This clearly implies that for any $\theta\in (0, \frac14)$ and $x_0\in B_{\frac14}$,
\beq{}\label{morrey_est2}
\|\psi\|_{M^{4,2}(B_\theta(x_0))}\le C\left(\epsilon_0^2+\theta^\frac12\right)\|\psi\|_{M^{4,2}(B_1)}.
\eeq
For any $\alpha\in (0,\frac13)$, first choose $\theta\in (0,\frac14)$ be such that $C\theta^{\frac12}
\le \theta^{\frac{\alpha}2}$ and then choose $\epsilon_0>0$ be such that
$C\epsilon_0^2\le \theta^{\frac{\alpha}2}$. Then we have
\beq{}\label{morrey_est3}
\|\psi\|_{M^{4,2}(B_\theta(x_0))}\le \theta^{\frac{\alpha}2}\|\psi\|_{M^{4,2}(B_1)},
\ \forall x_0\in B_{\frac14}.
\eeq
Iteration of (\ref{morrey_est3}) yields
\beq{}\label{morrey_est4}
\|\psi\|_{M^{4,2}(B_r(x_0))}\le Cr^{\frac{\alpha}2}\|\psi\|_{M^{4,2}(B_1)},
\ \forall x_0\in B_{\frac14}\ {\rm{ and } }\ 0\le r<\frac14.
\eeq
In particular, we have for any $0<\alpha<\frac13$,
\beq{}\label{morrey_est4}
r^{2(1-\alpha)-n}\int_{B_r(x_0)}|\psi|^4\le C\int_{B_1}|\psi|^4,
\ \forall x_0\in B_{\frac14} \ \hbox{ and } 0<r<\frac14.
\eeq
Thus $\psi\in M^{4,2(1-\alpha)}(B_\frac14)$ for any $\alpha\in (0,1)$.\\

\noindent{\it Step} 2 ($m=2$).
We want to repeat the above argument to show that
$\psi\in M^{\frac{4-4\alpha}{1-3\alpha}, 2(1-\alpha)}(B_{\frac1{16}})$.
In fact, since $(\eta_2|\psi|)^3\in M^{\frac43,2(1-\alpha)}(\mathbb R^n)$,
(\ref{adams}) implies that $f_2^i\in M^{\frac{4(1-\alpha)}{1-3\alpha}, 2(1-\alpha)}(\mathbb R^n)$,
and
\begin{eqnarray}\label{morrey_est5}
\left\|f_2^i\right\|_{M^{\frac{4(1-\alpha)}{1-3\alpha}, 2(1-\alpha)}(B_{\frac18})}
&\le& \left\|f_2^i\right\|_{M^{\frac{4(1-\alpha)}{1-3\alpha}, 2(1-\alpha)}(\mathbb R^n)}\nonumber\\
&\le& C\left\|\eta_2^3|\psi|^3\right\|_{M^{\frac43, 2(1-\alpha)}(\mathbb R^n)}\nonumber\\
&\le &C \left\|\psi\right\|_{M^{4,2(1-\alpha)}(B_{\frac14})}.
\end{eqnarray}
On the other hand, since $g_2^i$ is a harmonic function on $B_{\frac18}$, we have, by (\ref{morrey_est5}),
\begin{eqnarray}\label{morrey_est6}
&&\left\|g_2^i\right\|_{M^{\frac{4(1-\alpha)}{1-3\alpha}, 2(1-\alpha)}(B_{\frac1{16}})}\nonumber\\
&\le& C\left\|g_2^i\right\|_{M^{\frac{4(1-\alpha)}{1-3\alpha}, 2(1-\alpha)}(B_{\frac18})}\nonumber\\
&\le& C\left[\left\|f_2^i\right\|_{M^{\frac{4(1-\alpha)}{1-3\alpha}, 2(1-\alpha)}(B_{\frac18})}
+\left\|\psi^i\right\|_{M^{\frac{4(1-\alpha)}{1-3\alpha}, 2(1-\alpha)}(B_{\frac18})}\right]\nonumber\\
&\le & C\left\|\psi\right\|_{M^{\frac{4(1-\alpha)}{1-3\alpha}, 2(1-\alpha)}(B_{\frac18})}.
\end{eqnarray}
Putting (\ref{morrey_est5}) and (\ref{morrey_est6}) into (\ref{auxil2}) yields
that $\psi\in M^{\frac{4(1-\alpha)}{1-3\alpha}, 2(1-\alpha)}(B_{\frac1{16}})$ and
\beq{}\label{morrey_est7}
\left\|\psi\right\|_{M^{\frac{4(1-\alpha)}{1-3\alpha}, 2(1-\alpha)}(B_{\frac1{16}})}
\le C\left\|\psi\right\|_{M^{4,2(1-\alpha)}(B_{\frac14})} \le C\left\|\psi\right\|_{M^{4,2}(B_1)}.
\eeq
Since
$$\lim_{\alpha\uparrow\frac13}\frac{4(1-\alpha)}{1-3\alpha}=+\infty
\ {\rm{and}}\  \ M^{\frac{4(1-\alpha)}{1-3\alpha}, 2(1-\alpha)}(B_{\frac1{16}})
\subseteq L^{\frac{4(1-\alpha)}{1-3\alpha}}(B_{\frac1{16}}),$$
it follows from that $\psi\in L^p(B_{\frac1{16}})$ for any $p>4$,  and
\beq{}\label{lp_est}
\left\|\psi\right\|_{L^p(B_{\frac1{16}})}\le C(n,p)\|\psi\|_{M^{4,2}(B_1)}.
\eeq

Since
$|\partial\hspace{-0.25cm}/\hspace{0.125cm}\psi|
\le C|\psi|^3$,
$W^{1,p}$-estimate implies that $\psi\in W^{1,p}_{{\rm{loc}}}(B_{\frac1{16}},\mathbb C^{NL})$
for any $p>4$. Hence, by the Sobolev embedding theorem, $\psi\in C^\mu(B_{\frac1{16}},\mathbb C^{NL})$
for any $\mu\in (0,1)$. By the Schauder estimate, this yields $\psi\in C^{1,\mu}(B_{\frac1{16}},\mathbb C^{NL} )$.
Hence, by the bootstrap argument, we conclude $\psi\in C^\infty(B_{\frac1{16}},\mathbb C^{NL})$
and the estimate (\ref{prior-est}) holds.
\qed

\end{document}